\documentclass[11pt]{article}
\usepackage{amsfonts,amssymb,amsmath}

\newcommand{\cN}{{\cal N}}

\newcommand{\bbZ}{\mathbb Z}
\newcommand{\bbR}{\mathbb R}
\newcommand{\bbN}{\mathbb N}

\newcommand{\bbE}{\mathbb E}

\newcommand{\bbP}{\mathbb P}

\usepackage[english]{babel}
\newtheorem{theo}{Theorem}[section]
\newtheorem{pr}{Proposition}[section]
\newtheorem{cor}{Corollary}[section]
\newtheorem{lem}{Lemma}[section]
\newtheorem{defn}{Definition}[section]

\newcommand{\ra}{\mathop{\rightarrow }}

\newcommand{\Var}{\mathop{\rm Var}}
\newcommand{\Cov}{\mathop{\rm Cov}}

\newcommand{\epreuve}{\hspace{\fill}$\bigtriangleup$}

\begin{document}

\title{{\bf Central limit theorem for sampled sums of dependent random variables}}

\author{Nadine {\bf Guillotin-Plantard}\footnote{Universit\'e Claude Bernard - Lyon I, Institut Camille Jordan,
B\^atiment Braconnier, 43 avenue du 11 novembre 1918,  69622
Villeurbanne Cedex, France. {\it E-mail:}
nadine.guillotin@univ-lyon1.fr} and Cl\'ementine {\bf Prieur}\footnote{INSA Toulouse, Institut Math\'ematique de Toulouse, Equipe de Statistique et Probabilit\'es,
135 avenue de Rangueil, 31077 Toulouse Cedex 4, France.
{\it E-mail:} Clementine.Prieur@insa-toulouse.fr}}

\maketitle

\begin{abstract}
We prove a central limit theorem for linear triangular
arrays under weak dependence conditions. Our result is then applied
to the study of dependent random variables sampled by a
$\bbZ$-valued transient random walk. This extends the results
obtained by Guillotin-Plantard \& Schneider (2003). An application
to parametric estimation by random sampling is also provided.
\end{abstract}

\medskip


\vspace{7cm}



\newpage
\section{Introduction}
Let $\{\xi_i\}_{i \in \mathbb{Z}}$ be a sequence of centered, non
essentially constant and square integrable real valued random variables. Let
$\{a_{n,i} \, , \; -k_n \leq i \leq k_n \}$ be a triangular array of real numbers
such that for all $n \in \mathbb{N}$, $\sum_{i=-k_n}^{k_n}a_{n,i}^2 >0$. We
are interested in the behaviour of linear triangular arrays of the form
\begin{equation}\label{artri}
X_{n,i}=a_{n,i} \, \xi_i \, , \; n=0,1,\ldots ,\; i=-k_n, \ldots , k_n \, ,
\end{equation} where
$(k_n)_{n \geq 1}$ is a nondecreasing sequence of positive integers
satisfying $k_n \xrightarrow[n \rightarrow + \infty]{}+ \infty$. We
work under a weak dependence condition introduced in Dedecker {\it
et al.} (2007). We first prove a central limit theorem for linear
triangular arrays of type (\ref{artri}) (Theorem \ref{thmix0} of
Section \ref{clt}). Applying this result, we then prove a central
limit theorem for the partial sums of weakly dependent sequences sampled by a transient $\mathbb{Z}$-valued random walk
(Theorem \ref{thmix} of Section \ref{sampledclt}). This result
extends the results obtained by Guillotin-Plantard \& Schneider
(2003). Peligrad \& Utev (1997) derive a central limit theorem for
triangular arrays of type (\ref{artri}) under mixing conditions.
Unfortunately, mixing is a rather restrictive condition, and many
simple Markov chains are not mixing. For instance, Andrews (1984)
proved that if $(\varepsilon_i)_{i \geq 1}$ is independent and
identically distributed with marginal $\mathcal{B}\left(1/2\right)$,
then the stationary solution $(\xi_i)_{i \geq 0}$ of the equation
\begin{equation}\label{andr}
\xi_n=\frac12 (\xi_{n-1}+\varepsilon_n) \, , \; \xi_0 \textrm{ independent
  of } (\varepsilon_i)_{i \geq 1}
\end{equation}
is not $\alpha$-mixing in the sense of Rosenblatt (1956). We have indeed
  $\alpha(\sigma(\xi_0),\sigma(\xi_n))=1/4$ for any $n$. For any $y \in
  \mathbb{R}$, let $[y]$ denote the integer part of $y$. The chain
  satisfying (\ref{andr}) is the Markov chain associated to the dynamical
  system generated by the map $T(x)=2x-[2x]$ on the space $[0,1]$, equipped
  with the Lebesgue measure, and it is well known that such dynamical
  systems are not $\alpha$-mixing in the sense that
  $\alpha(\sigma(T),\sigma(T^n))$ does not tend to zero as $n$ tends to
  infinity. Withers (1981) proves triangular central limit theorems under
  a so-called $l$-mixing condition, which generalizes the classical notions
  of mixing (such as strong mixing, absolute regularity, uniform mixing
  introduced respectively by Rosenblatt (1956), Rozanov \& Volkonskii (1959)
  and Ibragimov (1962)). The idea of $l$-mixing requires the asymptotic
  decoupling of the 'past' and the 'future'. The dependence setting used in
  the present paper (introduced in Dedecker {\it et al.}, 2007) follows the same
  idea. In Section \ref{examples} we give lots of pertinent examples
  satisfying our dependence conditions. Coulon-Prieur \& Doukhan (2000)
  proves a triangular central limit theorem under a weaker dependence
  condition. However, they assume that the random variables $\xi_i$ are uniformly bounded. Their proof is a variation on Lindeberg-Rio's method developed
  by Rio (1996,1997). Also using a variation on Lindeberg-Rio's method, Bardet
  {\it et al.} prove a triangular central limit theorem, requiring moments of
  order $2+\delta$, $\delta>0$.
In Section \ref{def}, we introduce the dependence setting under
which we work in the sequel. Models for which we can compute bounds
for our dependence coefficients are presented in Section
\ref{examples}. At least, we give an application to parametric
estimation by random sampling in Section \ref{param}.
\section{Definitions}\label{def}
In this section, we recall the definition of the dependence coefficients
which we will use in the sequel. They have first been introduced in Dedecker
{\it et al.} (2007).

On the Euclidean space $\bbR^m$, we define the metric
$$d_1(x,y)= \sum_{i=1}^m |x_i - y_i|.$$
Let $\Lambda=\bigcup_{m\in\bbN^{*}}\Lambda_m$ where $\Lambda_m$ is the set of Lipschitz functions $f:\bbR^m\rightarrow \bbR$ with respect to the metric $d_1$.
If $f\in\Lambda_m$, we denote by $\mbox{\rm Lip}(f):=\sup_{x\neq y} \frac{|f(x)-f(y)|}{d_1(x,y)}$ the Lipschitz modulus of $f$. The set of functions $f\in \Lambda$ such that $\mbox{\rm Lip}(f)\leq 1$ is denoted by $\tilde\Lambda$.\\*

\begin{defn}\label{coeff}
 Let $\xi$ be a $\bbR^m$-valued random variable defined on a probability space
$(\Omega,{\cal A},\bbP)$, assumed to be square integrable.
For any $\sigma$-algebra ${\cal M}$ of ${\cal A}$, we define the $\theta_2$-dependence coefficient
\begin{equation}\label{tet2}
\theta_2({\cal M},\xi)=\sup \{\|\bbE(f (\xi)|{\cal M})-\bbE(f (\xi))\|_2 \, , \;
f\in \tilde\Lambda\} \, .
\end{equation}
\end{defn}

We now define the coefficient $\theta_{k,2}$ for a sequence of
$\sigma$-algebras and a sequence of $\mathbb{R}$-valued random variables.

\begin{defn}\label{coeffsuite}
Let $(\xi_i)_{i\in \bbZ}$ be a sequence of square integrable random variables
valued in $\mathbb{R}$. Let $({\cal M}_i)_{i\in\bbZ}$ be a sequence of
$\sigma$-algebras of ${\cal A}$. For any $k \in \mathbb{N}^*\cup \{\infty\}$
and $n \in \mathbb{N}$, we define
$$\theta_{k,2}(n)= \max_{1\leq l\leq k}
\frac 1 l \sup\{\theta_{2}({\cal M}_p, (\xi_{j_1}, \ldots ,\xi_{j_l})), p+n \leq j_1< \ldots <j_l\}$$ and
$$\theta_2(n)= \theta_{\infty,2} (n)= \sup_{k\in \mathbb{N}^*} \theta_{k,2}(n) \, .$$
\end{defn}
\begin{defn}
Let $(\xi_i)_{i\in \bbZ}$ be a sequence of square integrable random variables
valued in $\mathbb{R}$. Let $({\cal M}_i)_{i\in\bbZ}$ be a sequence of
$\sigma$-algebras of ${\cal A}$. The sequence $(\xi_i)_{i \in \bbZ}$ is said
to be $\theta_2$-weakly dependent with respect to $({\cal M}_i)_{i\in\bbZ}$
if $\theta_2(n) \xrightarrow[n \rightarrow + \infty]{} 0$.
\end{defn}

\noindent {\bf Remark:}\\
Replacing the
$\| \cdot \|_2$ norm in (\ref{tet2}) by the $\| \cdot \|_1$ norm, we get the
$\theta_1$ dependence coefficient first introduced by Doukhan \& Louhichi
(1999). This weaker coefficient is the one used in Coulon-Prieur \& Doukhan
(2000).

\section{Central limit theorem for triangular arrays of dependent random
  variables}\label{clt}
Let $\{ X_{n,i} , \ n \in \mathbb{N}, \ -k_n \leq i \leq k_n \}$ be a
triangular array of type (\ref{artri}).
We are interested in the asymptotic behaviour of the following sum
$$\Sigma_n=\sum_{i=-k_n}^{k_n} X_{n,i}=\sum_{i=-k_n}^{k_n}a_{n,i}\ \xi_{i} \, .$$
Let $({\cal M}_i)_{i \in \mathbb{Z}}$ be the sequence of
$\sigma$-algebras of ${\cal A}$ defined by
$$\mathcal{M}_i=\sigma\left(\xi_j \, , \; j \leq i\right)\, , \;
i\in \mathbb{Z} \, .$$ In the sequel, the dependence coefficients are
defined with respect to the sequence of $\sigma$-algebras $({\cal
M}_i)_{i \in \mathbb{Z}}$. We denote by $\sigma_n^2$ the variance of
$\Sigma_n$.
\begin{theo}\label{thmix0}
\hspace{1em}\\
Assume that the following conditions are satisfied~:
\begin{itemize}
\item[$(A_1)$]
\begin{enumerate}
\item[$(i)$] $\liminf_{n \rightarrow + \infty}
  \frac{\sigma_n^2}{\sum_{i=-k_n}^{k_n} a_{n,i}^2} >0$,
\item[$(ii)$] $
  \lim_{n \rightarrow + \infty} \sigma_n^{-1} \ \max_{-k_n \leq i \leq k_n}|a_{n,i}|= 0$.
\end{enumerate}
\item[$(A_2)$]
$\{\xi_i^2\}_{i \in \mathbb{Z}}$ is an uniformly integrable family.
\item[$(A_3)$] $\theta_2^{\xi}(\cdot)$ is bounded above by a
  non-negative function $g(\cdot)$ such that
\begin{itemize}
\item[] $x \mapsto x^{3/2} \ g(x)$ is non-increasing,
\item[] $\exists \ 0 < \varepsilon < 1 \ , \, \displaystyle\sum_{i=0}^{\infty}
  2^{\frac{3i}{2}}g(2^{i\varepsilon}) < \infty $.
\end{itemize}
\end{itemize}
Then, as $n$ tends to infinity, $ \frac{\Sigma_n}{\sigma_n}$
converges in distribution to ${\cN}(0,1)$.
\end{theo}

\vspace{1em}

\noindent {\bf Remark:}
\begin{itemize}
\item Theorem 2.2 (c) in Peligrad \& Utev (1997) yields a central limit theorem for
strongly mixing linear triangular arrays of type (\ref{artri}). They assume that
$\{|\xi_i|^{2+ \delta}\}$ is uniformly integrable for a certain $\delta
>0$.
Such an assumption is also required for Theorem 2.1 in Withers (1981) for
$l$-mixing arrays.
In Coulon-Prieur \& Doukhan (2000), the random variables $\xi_i$ are assumed to be
uniformly bounded.
\item The proof of Theorem 2.2 (c) in Peligrad \& Utev (1997) relies on a
  variation on Theorem 4.1 in Utev (1990) (see Theorem B in Peligrad \&
  Utev, 1997). The proof of Theorem \ref{thmix0}, which is postponed to the
  Appendix, also makes use of a variation on Theorem 4.1 in Utev (1990) (see
  also Utev, 1991).
\item If $\theta_2^{\xi}(n)=\mathcal{O}\left(n^{-a}\right)$ for some
  positive $a$, condition $(A_3)$ holds for $a>3/2$.
\end{itemize}
\section{Central limit theorem for the sum of dependent random variables sampled by a transient random walk}\label{sampledclt}
\subsection{The main result}
Let $(E,\mathcal{E},\mu)$ be a probability space, and $T  : E \mapsto E$ a bijective bimeasurable transformation preserving
the probability $\mu$.
We define the stationary sequence $(\xi_i)_{i
  \in \mathbb{Z}}=( T^i )_{i \in \mathbb{Z}}$ from $(E,\mu)$ to $E$.
Let $(X_i)_{i\geq
  1}$ be a sequence of independent and identically distributed random
variables defined on a probability space $(\Omega,\mathcal{A},\mathbb{P})$
with values in $\mathbb{Z}$ and
$$S_n=\sum_{i=1}^nX_i, \, n \geq 1, \quad \quad S_0 \equiv 0 \; .$$
For $f \in \mathbb{L}^1(\mu)$ and $\omega \in \Omega$, we are interested in
the sampled ergodic sum
$$\sum_{k=0}^{n-1}f \circ \xi_{S_k(\omega)} \, .$$
By applying Birkhoff's ergodic Theorem to the skew-product:
$$\begin{array}{rcl}
U & :  & \Omega \times E \rightarrow \Omega \times E\\
 & & (\omega,x) \mapsto (\sigma \omega,T^{\omega_1}x)
\end{array}$$
where $\sigma$ is the shift on the path space
$\Omega=\mathbb{Z}^{\mathbb{N}}$, we obtain that for every function $f \in
\mathbb{L}^1 (\mu)$, the sampled ergodic sum converges $\mathbb{P} \otimes
\mu$-almost surely. A natural question is to know if the random walk is
universally representative for $\mathbb{L}^p$, $p>1$ in the following sense:
there exists a subset $\Omega_0$ of $\Omega$ of probability one such that
for every $\omega \in \Omega_0$, for every dynamical system
$(E,\mathcal{E},\mu,T)$, for every $f \in \mathbb{L}^p$, $p>1$, the sampled
ergodic average converges $\mu$-almost surely. The answer can be found in
Lacey {\it et al.} (1994) if the $X_i$'s are square integrable: The random walk is
universally representative for $\mathbb{L}^p$, $p>1$ if and only if the
expectation of $X_1$ is not equal to $0$ which corresponds to the case where
the random walk is transient. In that case, it seems natural to study the
fluctuations of the sampled ergodic averages around the limit. From Lacey's
theorem (1991), for any $H \in (0,1)$, there exists some function $f \in
\mathbb{L}^2(\mathbb{P} \otimes \mu)$ such that the finite-dimensional
distributions of the process
$$\frac{1}{n^H} \sum_{k=0}^{[nt]-1} f \circ U^k (\omega,x)$$ converge to the
finite dimensional distributions of a self-similar process.
Unfortunately, this convergence on the product space does not imply
the convergence in distribution for a given path of the random walk.
A first answer to this question is given in Guillotin-Plantard \&
Schneider (2003) where the technique of martingale differences is
used. Let us recall that this method consists (under convenient
conditions) of decomposing the function $f$ as the sum of a function
$g$ generating a sequence of martingale differences and a cocycle
$h-h \circ T$. In the standard case, the central limit theorem for
the ergodic sum is deduced from central limit theorems for the sums
of martingale differences, the term corresponding to the cocycle
being negligeable in probability. In Guillotin-Plantard \& Schneider
(2003), only functions $f$ generating a sequence of martingale
differences are considered. In this section, where we prove a
central limit theorem for $\theta_2$-weakly dependent random
variables sampled by a transient random walk, this reasoning does
not hold anymore. We apply Theorem \ref{thmix0} of Section
\ref{clt}.

In the sequel, the random walk $(S_n)_{n\geq 0}$ is assumed to be transient. In particular,
for every $x \in \mathbb{Z}$, the Green function
$$ G(0,x)=\sum_{k=0}^{+\infty} \bbP(S_k=x)$$
is finite. For example, it is the case if the random variable $X_1$ is
assumed with finite absolute mean and nonzero mean. It is also possible to
choose the random variables $(X_i)_{i \geq 1}$ centered and for every $x \in
\mathbb{R}$,
$$\mathbb{P}(n^{-1/\alpha}S_n \leq x) \xrightarrow[n \rightarrow + \infty]{}
F_{\alpha}(x) \, ,$$
where $F_{\alpha}$ is the distribution function of a stable law with index
$\alpha \in (0,1)$. Stone (1966) has proved a local limit theorem for this
kind of random walks from which the transience can be deduced. The
expectation with respect to the measure $\mu$ (resp. with respect to
$\mathbb{P}$, $\mathbb{P} \otimes \mu$) will be denoted in the sequel by
$\mathbb{E}_{\mu}$ (resp. by $\mathbb{E}_{\mathbb{P}}$, $\mathbb{E}$).\\
For every function $f \in \mathbb{L}^2(\mu)$ such that $\mathbb{E}_{\mu}
(f)=0$, we define
$$\sigma^2(f)=2 \sum_{x \in \mathbb{Z}}G(0,x)\mathbb{E}_{\mu}(f f \circ T^x)-
\mathbb{E}_{\mu}(f^2) \, .$$ Let us now state our main result whose
proof is deferred to Subsection \ref{proofthmix}.

\begin{theo}\label{thmix}
\hspace{1em}\\
Let $f$ be a function in $\mathbb{L}^2(\mu)$ such that
$\mathbb{E}_{\mu}(f)=0$. Assume that $(f\circ T^x)_{x \in \mathbb{Z}}$
satisfies assumption $(A_3)$ of Theorem \ref{thmix0}.
Assume that $\sigma^2(f)$ is finite and positive.\\ Then, for
$\mathbb{P}$-almost every $\omega \in \Omega$, $$ \frac{1}{\sqrt{n}}
\sum_{k=0}^{n}f\circ T^{S_k(\omega)}  \xrightarrow[n \rightarrow +
\infty]{} {\cN}(0,\sigma^2(f)) \;\;\; \textrm{ in distribution.}$$
\end{theo}

\noindent {\bf Remark:}\\
\noindent {\bf 1.} In the particular case where $(f \circ T^x)_{x \in \mathbb{Z}}$ is a
sequence of martingale differences, we recognize Theorem 3.2 of
Guillotin-Plantard \& Schneider (2003). Indeed, assumptions are
satisfied using orthogonality of the $f \circ T^x$'s and then,
$\sigma^2(f)=(2G(0,0)-1)\mathbb{E}_{\mu}(f^2)$.\\
\noindent {\bf 2.} The stationarity assumption on the sequence $(\xi_i)_{i
  \in \mathbb{Z}}$ can be relaxed by a stationarity assumption of
order $2$, that is:
$$\forall \, i \in \mathbb{Z} \, , \; \Var \xi_i=\Var \xi_1 \; \; \textrm{and}$$
$$\forall \, i < j \, , \; \Cov (\xi_i,\xi_j)=\Cov (\xi_1,
\xi_{1+j-i}) \, .$$

\subsection{Computation of the variance}
The random walk $(S_n)_{n\geq 0}$ is defined as in the previous section.
The local time of the random walk is then defined for every $x\in\bbZ$ by
$$N_{n}(x)= \sum_{i=0}^{n} {\bf 1}_{\{S_{i}=x\}} \; .$$
The self-intersection local time is defined for every $x\in\bbZ$ by
$$\alpha(n,x)= \sum_{i,j=0}^{n} {\bf 1}_{\{S_{i}-S_{j}=x\}}$$
and can be rewritten using the definition of the local time as
$$\alpha(n,x)=\sum_{y\in\bbZ}N_{n}(y+x)\ N_{n}(y) \; .$$
Let $f$ be a function in $\mathbb{L}^2(\mu)$ such that
$\mathbb{E}_{\mu}(f)=0$. For every $\omega \in \Omega$,
$$\sum_{k=0}^{n}f \circ T^{S_k(\omega)}=\sum_{x \in \mathbb{Z}}N_n(x)(\omega)f \circ T^x
\,.$$ In order to apply results of Theorem \ref{thmix0}, we need to study,
for any fixed $\omega \in \Omega$, the asymptotic behaviour of the
variance of this sum, namely
$$\sigma_n^2(f)=\mathbb{E}_{\mu} \left( \left| \sum_{k=0}^n f \circ
T^{S_k(\omega)} \right|^2\right) \, .$$ The variable $\omega$ will
be omitted in the next calculations. We have
\begin{eqnarray*}
\sigma_{n}^2(f)&=&\bbE_{\mu}\left|\sum_{x\in\bbZ}N_{n}(x)f \circ T^{x} \right|^2\\
&=&\sum_{x,y\in\bbZ}N_{n}(x)N_{n}(y)\bbE_{\mu}(f \circ T^{x-y}\ f)\\
&=&\sum_{y,z\in\bbZ}N_{n}(y+z)N_{n}(y)\bbE_{\mu}(f \circ T^{z}\ f)\\
&=&\sum_{z\in\bbZ}\alpha(n,z)\bbE_{\mu}(f \circ T^{z}\ f)\, .
\end{eqnarray*}
We are now able to prove the following proposition:
\begin{pr}\label{pr2.1}
If $\sum_{x \in \mathbb{Z}} G(0,x) \mathbb{E}_{\mu}(f\ f \circ T^x)
< + \infty$, then $$\frac{\sigma_n^2(f)}{n} \xrightarrow[n
\rightarrow + \infty]{\mathbb{P}\textrm{-a.s.}}\sigma^2(f) \, .$$
\end{pr}
\noindent {\bf Proof of Proposition \ref{pr2.1}:}\\
Let us assume first that the function $f$ is positive. For every
$0\leq m<n$, we denote by $W_{m,n}$ the random variable
$$-\sum_{i,j=m}^{n}\sum_{x\in\bbZ}1_{\{S_{i}-S_{j}=x\}}\bbE_{\mu}(f\
f \circ T^{x}) \, .$$ Then, since $f$ is positive, for every $k,m,n$
such that $0\leq k<m<n$,
$$W_{k,n}\leq W_{k,m}+W_{m,n} \, ,$$
that is $(W_{m,n})_{m,n \geq 0}$ is a subadditive sequence. Then,
\begin{eqnarray*}
\bbE_{\mathbb{P}}(W_{0,n})&=&-\sum_{i,j=0}^{n}
\bbE_{\mathbb{P}}\Big(\sum_{x\in\bbZ}1_{\{S_{i}-S_{j}=x\}}\bbE_{\mu}(f\ f \circ T^{x})\Big)\\
&=&-\sum_{i,j=0}^{n}\sum_{x\in\bbZ}\bbP(S_{i}-S_{j}=x)\bbE_{\mu}(f\
f \circ T^{x}),\
\mbox{ by Fubini Theorem }\\
&=&-\Big((n+1)\bbE_{\mu}(f^2)+2\sum_{i=1}^{n}\sum_{j=0}^{i-1}\sum_{x \in \mathbb{Z}}\bbP(S_{i-j}=x)\bbE_{\mu}(f\ f
\circ T^{x})\Big)\\
&=&-\Big((n+1)\bbE_{\mu}(f^2)+2\sum_{i=1}^{n}\sum_{j=1}^{i} \sum_{x
\in \mathbb{Z}}\bbP(S_{j}=x)\bbE_{\mu}(f\ f \circ T^{x})\Big) \, .
\end{eqnarray*}
Now, using that
$$\lim_{i\ra +\infty}\sum_{j=1}^{i}\sum_{x \in \mathbb{Z}}\bbP(S_{j}=x)\bbE_{\mu}(f\ f \circ T^{x})
=\sum_{x\in\bbZ}G(0,x)\bbE_{\mu}(f\ f \circ T^{x})-\bbE_{\mu}(f^2)
,$$ we conclude that
$$\lim_{n \rightarrow + \infty}\frac{\bbE_{\mathbb{P}}(W_{0,n})}{n}=\bbE_{\mu}(f^2)
-2\sum_{x\in\bbZ}G(0,x)\bbE_{\mu}(f\ f \circ T^{x})<\infty \, .$$ So
the sequence $(W_{m,n})_{m,n\geq 0}$ satisfies all the conditions of
Theorem 5 in Kingman (1968). Hence
$$\frac{W_{0,n}}{n}\xrightarrow[n\ra + \infty]{\mathbb{P}\mbox{-a.s.}}\bbE_{\mu}(f^2)
-2\sum_{x\in\bbZ}G(0,x)\bbE_{\mu}(f\ f \circ T^{x}) \, .$$ By
remarking that $W_{0,n}=- \sigma_n^2(f)$, Proposition \ref{pr2.1}
follows for positive functions $f$. If the function $f$ is not
positive, we can decompose it as
$$f=f1_{\{f\geq 0\}}-(-f)1_{\{f<0\}}$$
and, for all $x \in \bbZ$, $f \circ T^x$ as
$$f \circ T^{x}=f \circ T^{x}1_{\{f \circ T^{x}\geq 0\}}-(-f \circ T^{x})1_{\{f \circ T^{x}<0\}} \, .$$
Then, a simple calculation yields
\begin{eqnarray*}
\bbE_{\mu}(f\ f \circ T^{x} )&=&\bbE_{\mu}(f\ f \circ T^{x}\
1_{\{f\geq 0;f \circ T^{x}\geq 0\}})
+\bbE_{\mu}((-f)(-f \circ T^{x})1_{\{f<0;f \circ T^{x}<0\}})\\
&-&\bbE_{\mu}(f\ (-f \circ T^{x})\ 1_{\{f\geq 0;f \circ T^{x}<0\}})-
\bbE_{\mu}((-f)\ (f \circ T^{x})\ 1_{\{f<0;f \circ T^{x}\geq 0\}})
\, .
\end{eqnarray*}
By applying the previous reasoning at each of the four positive
terms of the right-hand side, Proposition \ref{pr2.1} follows.
\epreuve

\noindent{\bf Remark:}\\
Let us consider the unsymmetric random walk on nearest neighbours
with $p>q$. Then, for $x \geq 0$,
$$G(0,x)= (p-q)^{-1} \, ,$$
and for $x\leq -1$,
$$G(0,x)= (p-q)^{-1}\left(\frac{p}{q}\right)^{x} \, .$$
A simple calculation gives
$$ \sigma^2(h-h \circ T)=2\sum_{x\in\bbZ} \left[ 2G(0,x) -G(0,x+1) -G(0,x-1)\right]
\bbE_{\mu}(h\ h \circ T^{x}) - 2 \bbE_{\mu}(h^2) + 2 \bbE_{\mu}(h\ h
\circ T)$$
$$ = - 2 \frac{p-1}{p}\bbE_{\mu}(h^2) + 2 \bbE_{\mu}(h\ h \circ T)-2 \frac{(p-q)}{pq} \sum_{x\geq 1}
\left(\frac{q}{p}\right)^{x} \bbE_{\mu}(h\ h \circ T^{x}) \, .$$

\subsection{\bf Proof of Theorem \ref{thmix}}\label{proofthmix}
Let us define $M_n=\displaystyle\max_{0\leq k\leq n} |S_k|$.
First note that
\begin{equation}\label{artriech}
\sum_{k=0}^{n} f \circ T^{S_k}=\sum_{|x|\leq M_n} N_n(x) f \circ T^{x} \, .$$
We want to apply Theorem \ref{thmix0} to the triangular array
$$\left\{X_{n,i}=\frac{N_n(i)}{\sqrt{n}} f \circ T^{i}, \ n \in \mathbb{N}, \
  -M_n \leq i \leq M_n \right\} \, .
\end{equation}
As $f$ belongs to $\mathbb{L}^2(\mu)$, the family $\left\{ \left( f \circ
  T^i \right)^2 \right\}_{i
  \in \mathbb{Z}}$ is uniformly integrable as it is stationary.
It remains to prove that assumption $(A_1)$ of Theorem
\ref{thmix0} is satisfied for the triangular array defined by
(\ref{artriech}).

\noindent {\bf Proof of $(A_1) (i)$:}\\
First, by Proposition 3.1. in Guillotin-Plantard \& Schneider
(2003),
$ \sum_{i=-M_n}^{M_n}a_{n,i}^2 = \frac{\alpha(n,0)}{n}$ converges
$\mathbb{P}$-almost surely to $2G(0,0)-1$ as $n$ goes to infinity.
Then, by Proposition \ref{pr2.1}, we know that $\sigma_n^2(f)/n$ converges to
$\sigma^2(f)$, which is assumed to be positive. Hence $(A_1) (i)$ is satisfied. \epreuve

\noindent {\bf Proof of $(A_1) (ii)$:}\\
Now, by Proposition 3.2. in Guillotin-Plantard \& Schneider
(2003), we know that for every $\rho>0$,
$$\max_{-M_n \leq i \leq M_n}|a_{n,i}|=\frac{1}{\sqrt{n}}\max_{i \in
   \bbZ}N_n(i)=o\left(n^{\rho-\frac{1}{2}}\right) \bbP-\textrm{almost
   surely.}$$So $\sqrt{\frac{n}{\sigma_n^2(f)}} \ \max_{-M_n \leq i \leq M_n}|a_{n,i}|$ tends to zero
$\bbP-$almost surely and assumption $(A_1) (ii)$ is satisfied. \epreuve

Hence Theorem \ref{thmix0} applied to
$\sum_{i=-M_n}^{M_n}a_{n,i}\ f \circ T^i$, Proposition \ref{pr2.1}
and Slutsky Lemma yield the result. \epreuve
\section{Examples}\label{examples}
In this section, we present examples for which we can compute upper bounds
for $\theta_{2}(n)$ for any $n \geq 1$. We refer to Chapter 3 in Dedecker \&
al (2007) and references therein for more details.
\subsection{Example 1: causal functions of stationary sequences}\label{ex1}
Let $(E,\mathcal{E},\mathbb{Q})$ be a probability space.
Let $(\varepsilon_i)_{i \in \mathbb{Z}}$ be a stationary sequence of random
variables with values in a measurable space $\mathcal{S}$. Assume that there
exists a real valued function $H$ defined on a subset of
$\mathcal{S}^{\mathbb{N}}$, such that $H(\varepsilon_0,\varepsilon_{-1},
\varepsilon_{-2}, \ldots,)$ is defined almost surely. The stationary
sequence $(\xi_n)_{n \in \mathbb{Z}}$ defined by $\xi_n=H(\varepsilon_n,
\varepsilon_{n-1}, \varepsilon_{n-2}, \ldots)$ is called a causal function
of $(\varepsilon_i)_{i \in \mathbb{Z}}$.

Assume that there exists a stationary sequence $({\varepsilon_i}')_{i \in
  \mathbb{Z}}$ distributed as $(\varepsilon_i)_{i \in \mathbb{Z}}$ and
  independent of $(\varepsilon_i)_{i \leq 0}$. Define $\xi_n^{*}=H({\varepsilon_n}',
{\varepsilon_{n-1}}', {\varepsilon_{n-2}}', \ldots)$. Clearly, $\xi_n^{*}$ is
  independent of $\mathcal{M}_0=\sigma(\xi_i \, , \; i \leq 0)$ and
  distributed as $\xi_n$. Let $(\delta_2(i))_{i >0}$ be a non increasing
  sequence such that
\begin{equation}\label{delta2}
\left\| \mathbb{E}\left(\left|\xi_i-\xi_i^{*}\right| \, | \,
  \mathcal{M}_0\right)\right\|_2 \leq \delta_2(i) \, .
\end{equation}
Then the
  coefficient $\theta_2$ of the sequence $(\xi_n)_{n \geq 0}$ satisfies
\begin{equation}\label{bound}
\theta_2(i) \leq \delta_2(i) \, .
\end{equation}
Let us consider the particular case where the sequence of innovations $(\varepsilon_i)_{i \in
  \mathbb{Z}}$ is absolutely regular in the sense of Volkonskii \& Rozanov
(1959). Then, according to Theorem 4.4.7 in Berbee (1979), if $E$ is
rich enough, there exists $(\varepsilon_i')_{i \in \mathbb{Z}}$ distributed
as $(\varepsilon_i)_{i \in \mathbb{Z}}$ and independent of
$(\varepsilon_i)_{i \leq 0}$ such that $$\mathbb{Q}(\varepsilon_i \neq
\varepsilon_i' \textrm { for some } i \geq k \, | \, \mathcal{F}_0)=\frac12
\left\| \mathbb{Q}_{\tilde{\varepsilon}_k|\mathcal{F}_0} -
\mathbb{Q}_{\tilde{\varepsilon}_k}\right\|_v ,$$
where $\tilde{\varepsilon}_k=(\varepsilon_k, \varepsilon_{k+1}, \ldots)$,
$\mathcal{F}_0=\sigma(\varepsilon_i \, , \, i \leq 0)$, and $\| \cdot \|_{v}$ is
  the variation norm. In particular if the
sequence $(\varepsilon_i)_{i \in \mathbb{Z}}$ is idependent and identically
distributed, it suffices to take $\varepsilon_i'=\varepsilon_i$ for $i>0$
and $\varepsilon_i'-\varepsilon_i''$ for $i \leq 0$, where
$(\varepsilon_i'')_{i \in \mathbb{Z}}$ is an independent copy of
$(\varepsilon_i)_{i \in \mathbb{Z}}$.

\noindent {\bf Application to causal linear processes:}\\
In that case, $\xi_n=\sum_{j \geq 0} a_j \varepsilon_{n-j}$, where $(a_j)_{j
  \geq 0}$ is a sequence of real numbers. We can choose
$$\delta_2(i) \geq \|\varepsilon_0-\varepsilon_0'\|_2 \sum_{j \geq i} |a_j| + \sum_{j=0}^{i-1} |a_j| \|\varepsilon_{i-j}- \varepsilon_{i-j}'\|_2 \, .
$$
From Proposition 2.3  in Merlev\`ede \& Peligrad (2002), we obtain
that
$$
  \delta_2(i) \leq  \|\varepsilon_0-\varepsilon_0'\|_2 \sum_{j \geq i} |a_j| + \sum_{j=0}^{i-1} |a_j| \,
\Big (2^2 \int_0^{\beta(\sigma(\varepsilon_k, k \leq 0), \sigma(\varepsilon_k, k \geq i-j))} Q^2_{\varepsilon_0}(u)\Big)^{1/2} du
 ,
$$
where $Q_{\varepsilon_0}$ is the generalized inverse of the tail function
$x \mapsto \mathbb{Q}(|\varepsilon_0|>x)$.

\subsection{Example 2: iterated random functions}
Let $(\xi_n)_{n \geq 0}$ be a real valued stationary Markov chain,
such that
$
\xi_n=F(\xi_{n-1}, \varepsilon_n)
$
for some measurable function $F$ and some independent and identically
distributed sequence $(\varepsilon_i)_{i > 0}$ independent of $\xi_0$.
Let $\xi^*_0$ be a random variable distributed as $\xi_0$ and independent
of
$(\xi_0, (\varepsilon_i)_{i >0})$.
Define
$
\xi^*_n=F(\xi^*_{n-1}, \varepsilon_n) \, .
$
The sequence $(\xi^*_n)_{n \geq 0}$ is distributed as $(\xi_n)_{n \geq  0}$
and independent of
$\xi_0$. Let ${\cal M}_i=\sigma(\xi_j , 0 \leq j \leq i)$.
As in Example 1,
define the sequence
$(\delta_2(i))_{i > 0}$ by (\ref{delta2}).
The coefficient  ${\theta}_2$ of
the sequence $(\xi_n)_{n \geq 0}$ satisfies the bound (\ref{bound}) of Example 1.

Let $\mu$ be the distribution of
$\xi_0$ and $(\xi_n^x)_{n \geq 0}$ be the chain starting from $\xi_0^x=x$.
With these notations, we can choose $\delta_2(i)$ such that
$$\delta_2(i) \geq \| \xi_i-\xi^*_i\|_2=\left(\int\int \|
|\xi_i^x-\xi_i^y\|_2^2 \mu(dx) \mu(dy) \right)^{1/2} \, .$$
For instance, if there exists a sequence
$(d_2(i))_{i \geq 0}$ of positive numbers
such that
$$
\|\xi_i^x-\xi_i^y\|_2 \leq d_2(i) |x-y|  ,
$$
then we can take
$ \delta_2(i) =d_2(i) \| \xi_0-\xi_0^* \|_2$.
For example, in the usual case where
$ \|F(x, \varepsilon_0)-F(y, \varepsilon_0)\|_2 \leq \kappa |x-y| $
for some $\kappa < 1$, we can
take $d_2(i)= \kappa^i$.

An important example is $\xi_n=f(\xi_{n-1}) + \varepsilon_n$
for some $\kappa$-Lipschitz function $f$. If $\xi_0$
has a moment of order $2$, then
$
\delta_2(i) \leq \kappa^i \|\xi_0-\xi_0^*\|_2 \, .
$

\subsection{Example 3: dynamical systems on $[0,1]$}
Let $I=[0,1]$, $T$ be a map from $I$ to $I$ and define $X_i=T^i$. If
$\mu$ is invariant by $T$, the sequence
$(X_i)_{i\geq0}$ of random variables from $(I,\mu)$ to $I$ is strictly stationary.

For any finite measure $\nu$ on $I$, we use the notations
$\nu(h)=\int_I h(x) \nu(dx)$.
For any finite signed measure $\nu$ on $I$, let $\|\nu\|=|\nu|(I)$ be the total variation
of $\nu$.
 Denote by $\|g\|_{1, \lambda}$ the ${\mathbb L}^1$-norm with respect to
the Lebesgue measure $\lambda$ on $I$. \medskip

\noindent{\bf Covariance inequalities.} In many interesting cases, one can prove that,
 for any $BV$ function $h$ and any $k$ in ${\mathbb L}^1(I, \mu)$,
\begin{equation}\label{sd1}
|\Cov(h(X_0),k(X_n))| \leq a_n  \|k(X_n)\|_1(\|h\|_{1, \lambda}+ \|dh\|) \, ,
\end{equation}
for some nonincreasing sequence $a_n$ tending to zero as $n$ tends to
infinity.

\medskip

\noindent{\bf Spectral gap.}
Define the operator $\mathcal{L}$ from ${\mathbb L}^1(I, \lambda)$ to ${\mathbb L}^1(I, \lambda)$ $via$ the equality
$$
\int_0^1 \mathcal{L}(h)(x) k(x) d \lambda (x) = \int_0^1h(x) (k \circ T)(x) d \lambda (x) \ \, \text{ where $h \in {\mathbb L}^1(I, \lambda)$ and
$k \in {\mathbb L}^\infty(I, \lambda)$}.
$$
The operator $\mathcal{L}$ is called the Perron-Frobenius operator of $T$.
In many interesting cases, the spectral analysis of ${\cal L}$  in the Banach space of $BV$-functions
equiped with the norm $\|h\|_v=\|dh\|+\|h\|_{1, \lambda}$ can be done by using the Theorem of
Ionescu-Tulcea and Marinescu (see Lasota and Yorke (1974) and Hofbauer and Keller (1982)).
Assume that $1$ is a simple eigenvalue of ${\cal L}$ and that the rest
of the spectrum is contained in a closed disk of radius strictly smaller than one.
Then there exists
a unique $T$-invariant absolutely continuous probability $\mu$ whose density $f_\mu$ is $BV$, and
\begin{equation}\label{sg1}
 {\cal L}^n(h)=\lambda(h) f_\mu + \Psi^n(h) \quad \text{with} \quad
\|\Psi^n(h)\|_v \leq K\rho^n \|h\|_v.
\end{equation}
for some $0\leq \rho <1$ and $K>0$.
Assume moreover that:
\begin{equation}\label{sg2}
\text{$I_*=\{f_\mu \neq 0\}$ is an interval, and there exists $\gamma>0$ such that $f_\mu > \gamma^{-1}$ on $I_*$.}
\end{equation}
Without loss of generality assume that $I_*=I$ (otherwise, take the restriction
to $I_*$ in what follows).
Define now the Markov kernel associated to $T$ by
\begin{equation}\label{Mk}
P(h)(x)=\frac{\mathcal{L}(f_\mu h)(x)}{f_\mu(x)} .
\end{equation}
It is easy to check (see for instance Barbour {\it et al.} (2000)) that $(X_0,X_1, \ldots , X_n)$
has the same distribution as $(Y_n,Y_{n-1},\ldots,Y_0)$ where $(Y_i)_{i \geq 0}$
is a stationary Markov chain with invariant distribution $\mu$ and transition kernel $P$.
Since $\|fg\|_\infty \leq \|fg\|_v \leq 2\|f\|_v \|g\|_v$, we infer that, taking $C=2K\gamma (\|df_\mu\|+1)$,
\begin{equation}\label{sg3}
 P^n(h)=\mu(h) + g_n \quad \text{with} \quad
\|g_n\|_\infty \leq C\rho^n \|h\|_v.
\end{equation}
This estimate implies   (\ref{sd1}) with $a_n=C\rho^n$ (see Dedecker \&
Prieur, 2005).

\medskip

\noindent {\bf Expanding maps:} Let  $([a_i, a_{i+1}[)_{1\leq i \leq N}$ be a finite partition of $[0,1[$.
We make the same assumptions on $T$ as in Collet {\it et al} (2002).
\begin{enumerate}
\item For each $1 \leq j \leq N$, the restriction $T_j$ of $T$ to $]a_j, a_{j+1}[$ is strictly monotonic
and can be extented to a function $\overline T_j$ belonging to $C^2([a_j, a_{j+1}])$.
\item Let $I_n$ be the set where $(T^n)'$ is defined.
There exists  $A>0$ and $s>1$ such that $ \inf _{x \in I_n} |(T^n)'(x)|>As^n$.
\item The map $T$ is topologically mixing: for any two nonempty open sets $U,V$, there exists $n_0\geq 1$
such that $T^{-n}(U)\cap V \neq \emptyset$ for all $n \geq n_0$.
\end{enumerate}
If $T$ satisfies 1.,2. and 3.,then (\ref{sg1}) holds. Assume furthermore
that (\ref{sg2}) holds (see Morita
(1994) for sufficient conditions). Then, arguing as in Example 4 in Section
7 of Dedecker \& Prieur (2005), we can prove that for the Markov chain
$(Y_i)_{i \geq 0}$ and the $\sigma$-algebras $\mathcal{M}_i=\sigma(Y_j \, ,
\, j \leq i)$, there exists a positive constant $C$ such that
$\theta_2(i) \leq C \rho^i $.

\medskip

\noindent {\bf Remark:}\\
In examples 2 and 3, the sequences are indexed by $\mathbb{N}$ and
not by $\mathbb{Z}$. However, using existence Theorem of Kolmogorov
(see Theorem 0.2.7 in Dacunha-Castelle \& Duflo, 1983), if $(X_i)_{i
\in \mathbb{N}}$ is a stationary process indexed by $\mathbb{N}$,
there exists a stationary sequence $(Y_i)_{i \in \mathbb{Z}}$
indexed by $\mathbb{Z}$ such that for any $k \leq l \in \mathbb{Z}$,
both marginals $(Y_k, \ldots , Y_l)$ and $(X_0, \ldots , X_{l-k})$
have the same distribution. Moreover, in examples 2 and 3, the
sequences are Markovian, hence $\theta_2^Y(n)=\theta_2^X(n)$ for any
$n \geq 1$. We then apply Theorem \ref{thmix} to the sequence
$(Y_i)_{i \in \mathbb{Z}}$. The limit variance can be rewritten as
$$\sigma^2(f)=2 \sum_{x \in \mathbb{Z}} G(0,x) \Cov (f(X_0),
f(X_{|x|})) - \Var (f(X_0)) \, .$$
\section{Application to parametric estimation by random
sampling}\label{param} We investigate in this section the problem of
parametric estimation by random sampling for second order stationary processes.
We assume that we observe a stationary process $(\xi_i)_{i \in
\mathbb{N}}$ at random times $S_n, \, n \geq 0$, where $(S_n)_{n
\geq 0}$ is a non negative increasing random walk satisfying the
assumptions of Section \ref{sampledclt}. In the case where the
marginal expectation of the process $(\xi_i)_{i \in \mathbb{N}}$, m,
is unknown, Deniau {\it et al.} estimate it using the sampled
empirical mean $\hat{m}_n=\frac1n \sum_{i=1}^{n} \xi_{S_i}$. They
measure the quality of this estimator by considering the following quadratic
criterion function:
$$a(S)=\lim_{n \rightarrow + \infty} (n \Var \hat{m}_n)\, .$$
In the case where $\left( \Cov(\xi_1, \xi_{n+1}) \right)_{n \in \mathbb{N}}$
is in $l^1$, we have
$$a(S)=\sum_{k=- \infty}^{+\infty} \Cov(\xi_{S_1},
\xi_{S_{|k|+1}}) <  \infty \, .$$

We then get Corollary \ref{corparam} below, which gives the asymptotic
behaviour of the estimate $\hat{m}_n$ after centering and normalization.
\begin{cor}\label{corparam}
Let us keep the assumptions of Section \ref{sampledclt} on the random walk
$(S_n)_{n \in \mathbb{N}}$ and on the process $(\xi_i)_{i \in
  \mathbb{N}}$. Assume moreover that $S_0=0$ and that $(S_{n+1}-S_n)_{n \in
  \mathbb{N}}$ takes its values in $\mathbb{N}^*$.
Then, for $\mathbb{P}$-almost every $\omega \in \Omega$,
$$\sqrt{n} \left( \hat{m}_n -m\right) \xrightarrow[n \rightarrow + \infty]{}
\mathcal{N}(0, a(S)) \, .$$
\end{cor}
\noindent {\bf Proof of Corollary \ref{corparam}:}\\
Corollary \ref{corparam} can be deduced from Theorem \ref{thmix} of Section
\ref{sampledclt} applied to $f=Id-m$. We have indeed $\sigma^2(f)=a(S)$.
\epreuve

\vspace{1em}

Let $\kappa \in \mathbb{R}_+$, $\kappa \geq 1$. A $\kappa$-optimal law for the step $S_{n+1}-S_n$ is a distribution which minimizes
$a(S)$ under the constraint $\mathbb{E} (S_{n+1}-S_n) \leq \kappa$. We say that
the  sampling
rate is less or equal to $1/\kappa$.
Deniau {\it et al.} (1988) give sufficient conditions for the existence of
a $\kappa$-optimal law of the step $S_{n+1}-S_n$. In the case of an A.R.(1)
model $\xi_n=\rho \xi_{n-1} + \varepsilon_n$, $|\rho | < 1$, with
$\varepsilon$ a white noise, we know from Taga (1965) that there exists a
unique $\kappa$-optimal law $L_0$ given by~:
\begin{itemize}
\item $L_0=\delta_1$ if $\rho < 0$ (the sampled process is the process
  itself),
\item $L_0=\left[1-\left(\kappa-[\kappa]\right)\right]\, \delta_{[\kappa]} +
  \left(\kappa-[\kappa]\right)\, \delta_{[\kappa]+1}$ if $\rho >0$
  \hspace{0.1em} ($[
  \cdot  ]$ denotes the integer part).
\end{itemize}

\section{Appendix}
This section is devoted to the proof of Theorem \ref{thmix0} of
Section \ref{clt}.

\noindent {\bf Proof of Theorem \ref{thmix0}:}\\
In order to prove Theorem \ref{thmix0}, we first use a classical
truncation argument. For any $M>0$, we define~:
\begin{equation*}
\varphi_{M}~:
\begin{cases}
\mathbb{R} \rightarrow \mathbb{R}\\
x \mapsto \varphi_{M}(x)=(x \wedge M) \vee (-M)
\end{cases}
\end{equation*}
and
\begin{equation*}
\varphi^M~:
\begin{cases}
\mathbb{R} \rightarrow \mathbb{R}\\
x \mapsto \varphi^M(x)=x-\varphi_{M}(x) \, .
\end{cases}
\end{equation*}
We now prove the following Lindeberg condition~:
\begin{equation}\label{lind}
\sigma_n^{-2}   \sum_{i=-k_n}^{k_n} \mathbb{E} \left( \left( \varphi^{\varepsilon
  \sigma_n}(X_{n,i})\right)^2 \right) \xrightarrow[n \rightarrow + \infty]{}
  0 \, .
\end{equation}
We have, for $n$ large enough,
$$\begin{array}{rcl}
\sigma_n^{-2} \sum_{i=-k_n}^{k_n} \mathbb{E} \left( \left( \varphi^{\varepsilon
  \sigma_n}(X_{n,i})\right)^2\right)   &
\leq & \sigma_n^{-2} \sum_{i=-k_n}^{k_n}
a_{n,i}^2 \mathbb{E} \left( \xi_i^2 {\bf 1}_{|\xi_i||a_{n,i}|> \varepsilon
    \sigma_n}\right)\\
& \leq &  \sigma_n^{-2} \sum_{i=-k_n}^{k_n}
a_{n,i}^2 \mathbb{E} \left( \xi_i^2 {\bf 1}_{|\xi_i| > \varepsilon
   \sigma_n /  \max_j |a_{n,j}|}\right) \, .
\end{array}$$
The last right hand term in the above inequalities is bounded by
\begin{equation}\label{ouf}
 \max_{-k_n \leq i \leq k_n} \left(\mathbb{E} \left( \xi_i^2 {\bf 1}_{|\xi_i| > \varepsilon
   \sigma_n /  \max_j |a_{n,j}|}\right)\right) \sigma_n^{-2}
\sum_{i=-k_n}^{k_n} a_{n,i}^2 \, ,
\end{equation}
which tends to zero as $n$ goes to infinity, by assumptions $(A_1)$ and
$(A_2)$.

By (\ref{ouf}) we find a sequence of positive numbers $(\varepsilon_n)_{n
  \geq 1}$ such that $\varepsilon_n \xrightarrow[n \rightarrow + \infty]{}
  0$, and
\begin{equation}\label{e1}
\max_{-k_n \leq i \leq k_n} \left(\mathbb{E} \left( \xi_i^2 {\bf 1}_{|\xi_i| > \varepsilon_n
   \sigma_n /  \max_j |a_{n,j}|}\right)\right) \sigma_n^{-2}
\sum_{i=-k_n}^{k_n} a_{n,i}^2 \xrightarrow[n \rightarrow + \infty]{} 0
\, .
\end{equation}
Let us now prove that (\ref{e1}) yields
\begin{equation}\label{lindbis}
\sigma_n^{-2} \Var \left( \sum_{i=-k_n}^{k_n}  \varphi^{\varepsilon_n
  \sigma_n}(X_{n,i})\right) \xrightarrow[n \rightarrow + \infty]{}
  0 \, .
\end{equation}
To prove (\ref{lindbis}), we need the following Lemma~:
\begin{lem}\label{maj1}
Assume that $(\eta_i)_{i \in \mathbb{Z}}$ is centered and satisfies conditions $(A_2)$ and
$(A_3)$ of Theorem \ref{thmix0}, then for any reals $-k_n \leq a \leq b
\leq k_n$,
$$\Var \left(\sum_{i=a}^b a_{n,i} \eta_i\right) \leq C \sum_{i=a}^b a_{n,i}^2
\, ,$$
with $C=\sup_{i \in \mathbb{Z}} \left(\mathbb{E}\eta_i^2\right) + 2 \sqrt{\sup_{i
    \in \mathbb{Z}} \left(\mathbb{E}\eta_i^2\right)} \ \sum_{l=1}^{\infty}
\theta_{1,2}(l)$.
\end{lem}

\vspace{1.5em}

Before proving Lemma \ref{maj1}, we finish the proof of (\ref{lindbis}).

For any fixed $n \geq 0$, and any $-k_n \leq i \leq k_n$ such that $a_{n,i}
\neq 0$, define~:
$$V_{n,i}=\varphi^{\varepsilon_n
    \sigma_n / |a_{n,i}|}(\xi_i)-\mathbb{E}\left( \varphi^{\varepsilon_n
    \sigma_n / |a_{n,i}|}(\xi_i) \right) \, .$$
If $a_{n,i}=0$, let $V_{n,i} =0$. As
    for any fixed $n \geq 0$, $-k_n \leq i \leq k_n$, the function
$$x \mapsto \varphi^{\varepsilon_n
    \sigma_n / |a_{n,i}|}(x)-\mathbb{E}\left( \varphi^{\varepsilon_n
    \sigma_n / |a_{n,i}|}(\xi_i) \right)$$
is $1$-Lipschitz, we have for all $l \geq 1$, for all $k \geq 1$,
$$\theta_{k,2}^{V_{\cdot ,n}}(l) \leq \theta_{k,2}^{\xi}(l) \, ,$$
where $V_{\cdot,n}=(V_{n,i})_{-k_n \leq i \leq k_n}$ and
$\xi=(\xi_i)_{i \in
  \mathbb{Z}}$.
Hence, for any fixed $n$, the sequence

$(V_{n,i})_{-k_n \leq i \leq k_n}$ satisfies
    assumptions $(A_2)$ and $(A_3)$ of Theorem \ref{thmix0}. Moreover,
as $$a_{n,i} \varphi^{\varepsilon_n
    \sigma_n / |a_{n,i}|}(\xi_i)=\varphi^{\varepsilon_n
    \sigma_n }(a_{n,i} \xi_i),$$ applying Lemma \ref{maj1} yields
\begin{equation}\label{majbis}
\sigma_n^{-2}\Var\Big(\sum_{i=-k_n}^{k_n} \varphi^{\varepsilon_n
    \sigma_n }(a_{n,i} \xi_i)\Big) \leq C_n \sigma_n^{-2}\sum_{i=-k_n}^{k_n} a_{n,i}^2
\, ,
\end{equation}
with $C_n=\sup_{-k_n \leq i \leq k_n} \left(\mathbb{E}V_{n,i}^2\right)
+ 2
    \sqrt{\sup_{-k_n \leq i
    \leq k_n} \left(\mathbb{E}V_{n,i}^2\right)} \ \sum_{l=1}^{\infty}
\theta_{1,2}^{\xi}(l)$. It remains to prove that the right hand term in
    (\ref{majbis}) converges to $0$ as
    $n$ goes to infinity.

We have, for $n$ large enough,
\begin{equation}\label{maj2}
\mathbb{E}V_{n,i}^2 \leq \mathbb{E} \left( \xi_i^2 {\bf 1}_{|\xi_i|
>
  \varepsilon_n  \sigma_n /  \max_j |a_{n,j}|} \right) \, .
\end{equation}
Hence we conclude with (\ref{e1}).

\vspace{1.5em}

\noindent {\bf Proof of Lemma \ref{maj1}:}
\begin{eqnarray*}
\Var\Big(\sum_{j=a}^{b} a_{n,j} \eta_j\Big)&=&\sum_{j=a}^{b}  a_{n,j}^2 \Var(\eta_j) +
\sum_{i=a}^{b}\sum_{j=a; j\neq i}^{b} a_{n,i}\ a_{n,j} \Cov(\eta_i,\eta_j)\\
&\leq & \sum_{j=a}^{b}  a_{n,j}^2 \Var(\eta_j) + \sum_{i=a}^{b}a_{n,i}^2\sum_{j=a; j\neq i}^{b}
|\Cov(\eta_i,\eta_j)|
\end{eqnarray*}
by remarking that $|a_{n,i}|\ |a_{n,j}| \leq \frac{1}{2} (a_{n,i}^2+a_{n,j}^2)$.\\*
Then for any $j>i$, using Cauchy-Schwarz inequality, we obtain that
\begin{eqnarray*}
|\Cov(\eta_i,\eta_j)|&=& \left|\bbE\left(\eta_i \; \bbE \left( \eta_j|{\cal M}_i\right)\right)\right| \\
&\leq & || \eta_i ||_{2} \ \left| \left| \bbE \left( \eta_j | {\cal M}_i
\right)  \right| \right|_{2}\\
&\leq & ||\eta_i ||_{2} \ \theta_{1,2}(j-i).
\end{eqnarray*}
As $(\eta_i)_{i \in \mathbb{Z}}$ is centered, and as $(\eta_i^2)_{i \in \mathbb{Z}}$ is uniformly integrable, we deduce that
\begin{equation*}
\Var\Big(\sum_{j=a}^{b} a_{n,j} \eta_j\Big) \leq
C \sum_{j=a}^{b}  a_{n,j}^2 \, ,
\end{equation*}
with $C=\sup_{i \in \mathbb{Z}} \left(\mathbb{E}\eta_i^2\right) + 2 \sqrt{\sup_{i
    \in \mathbb{Z}} \left(\mathbb{E}\eta_i^2\right)} \ \sum_{l=1}^{\infty}
\theta_{1,2}(l)$ which is finite from assumptions $(A_2)$ and
$(A_3)$. \epreuve

\vspace{1.5em}

Define $Z_{n,i}$ by
$$\frac{\varphi_{\varepsilon_n \sigma_n}(X_{n,i})- \mathbb{E}\left(\varphi_{\varepsilon_n \sigma_n}(X_{n,i})\right)}{ \sqrt{\Var \left(\sum_{i=-k_n}^{k_n}\varphi_{\varepsilon_n
  \sigma_n}(X_{n,i})\right)}}\equiv \frac{\varphi_{\varepsilon_n \sigma_n}(X_{n,i})- \mathbb{E}\left(\varphi_{\varepsilon_n \sigma_n}(X_{n,i})\right)}{\sigma'_n}\, .$$
By (\ref{lindbis}) we conclude that, to prove Theorem \ref{thmix0}, it is
enough to prove it for the truncated sequence $\left(Z_{n,i}\right)_{n\geq
  0, -k_n \leq i \leq k_n}$, that is to show that
  \begin{equation}\label{tlctronc}
  \sum_{i=-k_n}^{k_n}Z_{n,i}
\xrightarrow[ n \rightarrow +
  \infty]{\mathcal{D}}\mathcal{N}(0,1) \, .
  \end{equation}

The proof is now a variation on the proof of Theorem 4.1 in
Utev (1990). Let
$$d_t(X,Y)=\left|\mathbb{E}e^{itX} - \mathbb{E} e^{itY} \right| \, .$$
To prove Theorem \ref{thmix0}, it is enough to prove that for all $t$,
$$d_t \left( \sum_{i=-k_n}^{k_n} Z_{n,i}, \eta \right) \xrightarrow[n
  \rightarrow + \infty]{} 0 \, ,$$
with $\eta$ the standard normal distribution. We first need some simple
  properties of the distance $d_t$. Let $X,X_1,X_2,Y_1,Y_2$ be random
  variables with zero means and finite second moments. We assume that the
  random variables $Y_1,Y_2$ are independent and that the distribution of
  $X_j$ coincides with that of $Y_j$, $j=1,2$.
We define $A_t(X)=d_t\left(X , \eta \sqrt{\mathbb{E}X^2}\right)$.
We have then the following inequalities~:
\begin{lem}[Lemma 4.3 in Utev, 1990]\label{ateate}
$$A_t(X)\leq \frac{2}{3}|t|^3 \mathbb{E}|X|^3,$$
$$A_t(Y_1+Y_2) \leq A_t(Y_1)+A_t(Y_2),$$
$$d_t(X_1+X_2,X_1) \leq \frac{t^2}{2} \left( \mathbb{E} X_2^2 +
  (\mathbb{E}X_1^2 \mathbb{E} X_2^2)^{1/2}\right),$$
$$d_t(\eta a, \eta b) \leq \frac{t^2}{2} |a^2-b^2|.$$
\end{lem}

We next need the following lemma~:
\begin{lem}\label{ut}
Let $0 < \varepsilon < 1$.
There exists some positive constant $C(\varepsilon)$ such that for all $a \in
\mathbb{Z}$, for all $v \in \mathbb{N}^*$, $A_t\left(\sum_{i=a+1}^{a+v}
Z_{n,i} \right)$ is bounded by
$$ C(\varepsilon)\left( |t|^3
  h^{2/\varepsilon} \sum_{i=a+1}^{a+v} \mathbb{E}\left(|Z_{n,i}|^3\right) +
  t^2 \left(h^{\frac{\varepsilon -1}{2}} + \sum_{j~: 2^j \geq h^{1/\varepsilon}}
  2^{\frac{3 j}{2}}g(2^{j\varepsilon})\right) {\sigma'}_n^{-2}
  \sum_{i=a+1}^{a+v}a_{n,i}^2\right),$$
where $h$ is an arbitrary positive natural number and with $g$ introduced in
  Assumption $(A_3)$ of Theorem \ref{thmix0}.
\end{lem}

\vspace{1.5em}

Before proving Lemma \ref{ut}, we achieve the proof of Theorem \ref{thmix0}.
By Lemma \ref{ut}, we have
$$
d_t \left( \sum_{i=-k_n}^{k_n} Z_{n,i} , \eta \right)=A_t\left(\sum_{i=-k_n}^{k_n}
Z_{n,i}\right)
 \leq  C(t,\varepsilon) \left( h^{2 / \varepsilon }\sum_{i=-k_n}^{k_n} \mathbb{E}(|Z_{n,i}|^3)
+ \delta(h)\right) ,$$ 
with $\delta(h)=\Big(h^{\frac{\varepsilon
-1}{2}}+\sum_{j~: 2^j \geq h^{\frac{1}{\varepsilon}}}
  2^{\frac{3 j}{2}}g(2^{j\varepsilon})\Big) {\sigma'}_n^{-2}
  \sum_{i=-k_n}^{k_n}a_{n,i}^2$.

\vspace{1em}

From (\ref{lindbis}) we deduce that $\sigma_n/{\sigma'_n} \xrightarrow[ n
   \rightarrow + \infty]{} 1 $.

Hence using assumptions $(A_1)(i)$ and
 $(A_3)$ we get $\delta(h)\xrightarrow[ h
   \rightarrow + \infty]{} 0 $.

On the other hand we have
\begin{equation}\label{hip}
\sum_{i=-k_n}^{k_n} \mathbb{E}(|Z_{n,i}|^3) \leq 2 \ \frac{\varepsilon_n
\sigma_n}{\sigma'_n} \sum_{i=-k_n}^{k_n}\Var(Z_{n,i}).
\end{equation}
For any fixed $n \geq 0$, and any $-k_n \leq i \leq k_n$ such that $a_{n,i}
\neq 0$, define~:
$$W_{n,i}=\varphi_{\varepsilon_n
    \sigma_n / |a_{n,i}|}(\xi_i)-\mathbb{E}\left( \varphi_{\varepsilon_n
    \sigma_n / |a_{n,i}|}(\xi_i) \right) \, .$$
If $a_{n,i}=0$, let $W_{n,i} =0$. As
    for any fixed $n \geq 0$, $-k_n \leq i \leq k_n$, the function
$$x \mapsto \varphi_{\varepsilon_n
    \sigma_n / |a_{n,i}|}(x)-\mathbb{E}\left( \varphi_{\varepsilon_n
    \sigma_n / |a_{n,i}|}(\xi_i) \right)$$
is $1$-Lipschitz, we have for all $l \geq 1$, for all $k \geq 1$,
$$\theta_{k,2}^{W_{\cdot ,n}}(l) \leq \theta_{k,2}^{\xi}(l) \, ,$$
where $W_{\cdot,n}=(W_{n,i})_{-k_n \leq i \leq k_n}$ and
$\xi=(\xi_i)_{i \in
  \mathbb{Z}}$.
Hence, arguing as for the proof of Lemma \ref{maj1}, we prove that for any fixed $n$, the sequence
$(W_{n,i})_{-k_n \leq i \leq k_n}$ satisfies
    assumptions $(A_2)$ and $(A_3)$ of Theorem \ref{thmix0}.
Therefore, for any
reals $-k_n \leq a \leq b \leq k_n$,
\begin{equation}\label{maj3}
\Var \left( \sum_{i=a}^b Z_{n,i} \right) \leq C {\sigma'}_n^{-2} \sum_{i=a}^b a_{n,i}^2
\end{equation}
with $C =\sup_{i \in \mathbb{Z}} \left(2\mathbb{E}\xi_i^2\right) + 2 \sqrt{\sup_{i
    \in \mathbb{Z}} \left(2 \mathbb{E}\xi_i^2\right)} \ \sum_{l=1}^{\infty}
\theta_{1,2}(l)$ which is finite from assumptions $(A_2)$ and
$(A_3)$.
Applying (\ref{maj3}) with $a=b$, we get that the right hand term of (\ref{hip}) is bounded by
$2 \ \varepsilon_n
\frac{\sigma_n}{\sigma'_n} \sum_{i=-k_n}^{k_n}\frac{C
    a_{n,i}^2}{{\sigma'}_n^2}$, which tends to zero as $n$ tends to
    infinity, using assumption $(A_1)(i)$ and the fact that $\varepsilon_n
    \xrightarrow[n \rightarrow + \infty]{}0$.

Consequently
$$\inf_{h \geq 1} \left(h^{2 / \varepsilon }\sum_{i=-k_n}^{k_n} \mathbb{E}(|Z_{n,i}|^3)
+ \delta(h)\right) \xrightarrow[n \rightarrow + \infty]{}0.$$
It achieves the proof of Theorem \ref{thmix0}.
\epreuve

\vspace{1.5em}

\noindent {\bf Proof of Lemma \ref{ut}:} Let $h \in \mathbb{N}^*$.
Let $0 < \varepsilon < 1$.  In the following, $C$,
$C(\varepsilon)$ denote constants which may vary from line to line.
Let $\kappa_{\varepsilon}$ be a positive constant greater than $1$
which will be precised further. Let $v < \kappa_{\varepsilon}\
h^{\frac{1}{\varepsilon}}$. We have
\begin{equation}\label{petit}
A_t \left(\sum_{i=a+1}^{a+v} Z_{n,i}\right)  \leq \frac{2}{3} |t|^3
\mathbb{E}\left|\sum_{i=a+1}^{a+v} Z_{n,i}\right|^3 \leq \frac{2}{3} \
\kappa_{\varepsilon}^2 \ |t|^3 \
h^{2/\varepsilon}\sum_{i=a+1}^{a+v} \mathbb{E}(|Z_{n,i}|^3)
\end{equation}
since $|x|^3$ is a convex function.

Let now $v \geq \kappa_{\varepsilon}\
h^{\frac{1}{\varepsilon}}$. Without loss of generality, assume that $a=0$. Let 
$\delta_{\varepsilon}=(1-\varepsilon^2 +2\varepsilon)/2$. Define then
$$m=[v^{\varepsilon}] , \ B= \left\{u \in \mathbb{N}~: \ 2^{-1} (v
-[v^{\delta_{\varepsilon}}]) \leq um \leq 2^{-1}v\right\},$$
$$A=\left\{ u \in \mathbb{N}~:0 \leq u \leq v, \
\sum_{i=um+1}^{(u+1)m}a_{n,i}^2 \leq (m/v)^{\varepsilon} \sum_{i=1}^v
a_{n,i}^2 \right\}.$$
Following Utev (1991) we prove that, for $0 < \varepsilon < 1$, $A \cap B$
is not wide for $v$ greater than $\kappa_{\varepsilon}$. We have indeed
$$|A \cap B| =|B|-|\overline{A}\cap B|\geq |B|-|\overline{A}|\geq
\frac{v^{(1-\varepsilon^2)/2}}{2}\Big(1-4v^{-(1-\varepsilon)^2/2}\Big)-\frac{3}{2} \ , $$
where $\overline{A}$ denotes the complementary of the set $A$.
We can find $\kappa_{\varepsilon}$ large enough so that $|A \cap B|$ be positive.

Let $u \in A \cap B$. We start from
the following simple identity
\begin{eqnarray}
\nonumber Q & \equiv & \sum_{i=1}^v
Z_{n,i}\\
\nonumber & = &
\sum_{i=1}^{um}Z_{n,i}+\sum_{i=um+1}^{(u+1)m}Z_{n,i}+\sum_{i=(u+1)m+1}^v
Z_{n,i}\\
\label{4.7} & \equiv & Q_1 + Q_2 + Q_3.\end{eqnarray}
By Lemma \ref{ateate},
\begin{equation}\label{4.8}
d_t(Q,Q_1+Q_3)=d_t(Q,Q-Q_2)\leq \frac{t^2}{2} \left(\mathbb{E} Q_2^2 +
(\mathbb{E}Q_2^2 \mathbb{E} Q^2)^{1/2}\right).
\end{equation}
Using (\ref{4.8}) and (\ref{maj3}), we get
\begin{equation}\label{4.9}
d_t(Q,Q_1+Q_3) \leq C t^2 v^{\frac{(\varepsilon-1)\varepsilon}{2}} {\sigma'}_n^{-2} \sum_{i=1}^v a_{n,i}^2.
\end{equation}
Now, given the random variables $Q_1$ and $Q_3$, we define two independent
random variables $g_1$ and $g_3$ such that the distribution of $g_i$
coincides with that of $Q_i$, $i=1,3$.
We have
$$\begin{array}{rcl}
d_t(Q_1+Q_3,g_1+g_3) & = &
\left|\mathbb{E}(e^{itQ_1}-1)(e^{itQ_3}-1)-\mathbb{E}(e^{itQ_1}-1)\mathbb{E}(e^{itQ_3}-1)\right|\\\\
& \leq  &  \left\|\mathbb{E}(e^{itQ_1}-1)\right\|_2 \left\|
\mathbb{E}\left(e^{itQ_3}-1-\mathbb{E}(e^{itQ_3}-1) \ | \,
\mathcal{M}_{um}\right)\right\|_2\\\\
& \leq & 2 |t| \left\|\sum_{i=1}^{um}Z_{n,i}\right\|_2 v \ |t| \ {\sigma'}_n^{-1}\left(
        \sum_{i=(u+1)m+1}^v
|a_{n,i}|\right) \theta_2^{\xi}(m+1)\\\
& \leq & C \ t^2 \ v^{3/2} \ {\sigma'}_n^{-2} \left(\sum_{i=1}^v
a_{n,i}^2\right) g(v^{\varepsilon}),
\end{array}$$
by (\ref{maj3}), Definition \ref{coeffsuite} and Assumption $(A_3)$ of
Theorem \ref{thmix0}.
Hence
\begin{equation}\label{4.10}
d_t(Q_1+Q_3,g_1+g_3) \leq  C t^2 f(v) {\sigma'}_n^{-2}
\sum_{i=1}^v a_{n,i}^2,
\end{equation}
where $f(v)=v^{3/2} \ g(v^{\varepsilon})$ is non-increasing by assumption
$(A_3)$ of Theorem \ref{thmix0}.

We also have by Lemma \ref{ateate}
\begin{equation}\label{4.11}
A_t(g_1+g_3)\leq A_t(g_1)+A_t(g_3).
\end{equation}
Finally, still by Lemma \ref{ateate}, and using Definition \ref{coeffsuite}, we have
\begin{eqnarray}
\nonumber d_t \left( \eta \sqrt{\mathbb{E}(Q^2)}, \eta
\sqrt{\mathbb{E}\left( (g_1+g_3)^2 \right) } \right) & \leq & \frac{t^2}{2}
\left| \mathbb{E} (Q^2) - \mathbb{E} \left( (g_1+g_3)^2 \right)  \right|\\
 \nonumber & \leq & \frac{t^2}{2} \left| \mathbb{E} (Q_2^2) + 2 \mathbb{E}
(Q_1Q_2) + 2 \mathbb{E}(Q_2Q_3) +2 \mathbb{E}(Q_1Q_3)\right|\\
\label{4.12} & \leq & C t^2 \left(  v^{\frac{(\varepsilon -1)\varepsilon}{2}} + f(v)\right) {\sigma'}_n^{-2}
\sum_{i=1}^v a_{n,i}^2.
\end{eqnarray}
Combining (\ref{4.9})-(\ref{4.12}), we get the following recurrent
inequality~:
$$\begin{array}{rcl}
A_t\left(\sum_{i=1}^v Z_{n,i}\right) & \leq & A_t\left(\sum_{i=1}^{um}
Z_{n,i}\right) + A_t\left(\sum_{i=(u+1)m+1}^{v} Z_{n,i}\right)\\
& + & C t^2 \left(  v^{\frac{(\varepsilon -1)\varepsilon}{2}} + f(v)\right) {\sigma'}_n^{-2}
\sum_{i=1}^v a_{n,i}^2
\end{array}$$
for $v \geq \kappa_{\varepsilon} \
h^{\frac{1}{\varepsilon}}\geq \kappa_{\varepsilon}$.

We then need the following Lemma, which is a variation on Lemma 1.2. in Utev
(1991).

\begin{lem}\label{ututut}
For every $\varepsilon\in\ ]0,1[$, denote $\delta_{\varepsilon}=(1-\varepsilon^2+2\varepsilon)/2$.
Let a non-decreasing sequence of non-negative numbers $a(n)$ be
specified, such that there exist non-increasing sequences of
non-negative numbers $\varepsilon(k)$, $\gamma(k)$ and a sequence of
naturals $T(k)$, satisfying conditions
$$T(k) \leq 2^{-1}(k+[k^{\delta_{\varepsilon}}]),$$
$$a(k) \leq \max_{k_0 \leq s \leq k}\left(a(T(s))+\gamma(s)\right)$$
for any $k \geq k_0$ with an arbitrary $k_0 \in \mathbb{N}^*$. Then
$$a(n) \leq a(n_0)+2 \sum_{k_0 \leq 2^j \leq n}\gamma(2^j) \ ,$$
for any $n \geq k_0$, where one can take $n_0=2^c$ with $c>\frac{2-\delta_{\varepsilon}}{1-\delta_{\varepsilon}}$.
\end{lem}

\vspace{1.5em}

\noindent {\bf Proof of Lemma \ref{ututut}:} The proof follows essentially
the same lines as the proof of Lemma 1.2. in Utev (1991) and therefore is
omitted here. \epreuve

\vspace{1.5em}

We now apply Lemma \ref{ututut} above with
\begin{itemize}
\item[$\star$] $k_0=\kappa_{\varepsilon} \ h^{\frac{1}{\varepsilon}}$,
\item[$\star$] for $k \ge k_0$, $T(k)=\max \left\{ u_km_k, k-u_km_k-m_k \right\}$ where $u_k$
  and $m_k$ are defined from $k$ as $u$ and $m$ from $v$ (see the proof of
  $A \cap B $ not wide),
\item[$\star$] $c<\frac{\ln(\kappa_{\varepsilon})}{\ln(2)}\ \ \ $ (we may
  need to enlarge $\kappa_{\varepsilon}$),
\item[$\star$] for $s \geq k_0$, $\gamma(s)=C  \ t^2 \ \left( s^{\frac{\varepsilon ( \varepsilon
-1)}{2}}+f(s) \right)$,
\item[$\star$] for $s \geq k_0$, $a(s)=\displaystyle\sup_{l\in\mathbb{Z}} \ \max_{k_0 \leq i \leq s} \
\frac{A_t\left(\sum_{j=l+1}^{l+i}Z_{n,j}\right)}{{\sigma'}_n^{-2}\sum_{j=l+1}^{l+i}a_{n,j}^2}
$.
\end{itemize}
Applying Lemma \ref{ututut} yields the statement of Lemma \ref{ut}. 
\epreuve


\begin{thebibliography}{99}
\small
\bibitem{A} D. W. K. Andrews (1984). Nonstrong mixing autoregressive processes.
{\it J. Appl. Probab.} 21, p. 930-934.
\bibitem{BGR} A. D. Barbour, R. M. Gerrard and G. Reinert (2000). Iterates of expanding maps.
{\it Probab. Theory Relat. Fields} 116, p. 151-180.
\bibitem{Bar} J. M. Bardet, P. Doukhan, G. Lang and N. Ragache
  (2007). Dependent Linderberg central limit theorem and some applications. To
  appear in {\it ESAIM P\&S}.
\bibitem{Bee} H. C. P. Berbee (1979). Random walks with stationary increments and renewal theory. {\it Math. Cent. Tracts.} Amsterdam.

\bibitem {CMS} P. Collet, S. Martinez and B. Schmitt (2002). Exponential inequalities for
dynamical measures of expanding maps of the interval. {\it
  Probab. Theory. Relat. Fields} 123, p. 301-322.
\bibitem {CD} C. Coulon-Prieur and P. Doukhan (2000). A triangular CLT
 for weakly dependent sequences.
{\it  Statist.   Probab. Lett.} 47, p. 61-68.
\bibitem{DaDu} D. Dacunha-Castelle and M. Duflo. Probl\`emes \`a temps
  mobile.  Deuxi\`eme \'edition, Masson (1993).
\bibitem{DP05} J. Dedecker and C. Prieur (2005). New dependence
coefficients. Examples and applications to statistics. {\it Probab. Theory
  Relat. Fields} 132, p. 203-236.
\bibitem{Dedal} J. Dedecker, P. Doukhan, G. Lang, J. R. Leon, S. Louhichi
  and C. Prieur. Weak dependence: With Examples and Applications.
 Lect. notes in Stat. 190.
Springer, XIV (2007).
\bibitem{DOV} C. Deniau, G. Oppenheim and M. C. Viano (1988). Estimation de param\`etre par \'echantillonnage al\'eatoire. (French. English summary) [Random sampling and parametric estimation].
{\it C. R. Acad. Sci. Paris S\'er. I Math.} 306, 13, p. 565--568.
\bibitem {DL} P. Doukhan and S. Louhichi (1999).  A new weak dependence condition and applications to moment inequalities. {\it Stochastic Process. Appl.} 84, p. 313-342.
\bibitem{GPS} N. Guillotin-Plantard and D. Schneider (2003). Limit theorems
  for sampled dynamical systems. {\it Stochastic and Dynamics} 3, 4,
  p. 477-497.
\bibitem{HK82} F. Hofbauer and G. Keller (1982). Ergodic properties of
  invariant measures for piecewise monotonic transformations. {\it Math. Z.}
  180, p. 119-140.
\bibitem{Ib} I. A. Ibragimov (1962). Some limit theorems for stationary
  processes. {\it Theory Probab. Appl.} 7, p. 349-382.
\bibitem{Kin} J. F. C. Kingman (1968). The ergodic theory of subadditive
  stochastic processes. {\it J. R. Statist. Soc.} B30, p. 499-510.
\bibitem{Lac} M. Lacey (1991). On weak convergence in dynamical systems to
  self-similar processes with spectral representation. {\it
  Trans. Amer. Math. Soc.} 328, p. 767-778.
\bibitem{Lacal} M. Lacey, K. Petersen, D. Rudolph and M. Wierdl
  (1994). Random ergodic theorems with universally representative
  sequences. {\it Ann. Inst. H. Poincar\'e Probab. Statist.} 30, p. 353-395.
\bibitem{LaYo} A. Lasota and J. A. Yorke (1974). On the existence of
  invariant measures for piecewise monotonic transformations. {\it
  Trans. Amer. Math. Soc.} 186, p. 481-488.
\bibitem {MP} F. Merlev\`ede and  M. Peligrad (2002).
 On the coupling of dependent random variables and applications.
{\it Empirical process techniques for dependent data.}, p. 171-193. Birkh\"auser.
\bibitem{Mo} T. Morita (1994). Local limit theorem and distribution of
  periodic orbits of Lasota-Yorke transformations with infinite Markov
  partition. {\it J. Math. Soc. Japan} 46, p. 309-343.
\bibitem{PU} M. Peligrad and S. Utev (1997). Central limit theorem for
  linear processes. {\it Ann. Probab.} 25, 1, p. 443-456.
\bibitem{Ri1} E. Rio (1996). Sur le th\'eor\`eme de Berry-Esseen pour les
  suites faiblement d\'ependantes. {\it Probab. Theory Relat. Fields} 104,
  p. 255-282.
\bibitem{Ri2} E. Rio (1997). About the Lindeberg method for strongly mixing
  sequences. {\it ESAIM P\&S} 1, p. 35-61.
\bibitem{Ro} M. Rosenblatt (1956). A central limit theorem and a strong
  mixing condition. {\it Proc. Nat. Acad. Sci. U. S. A. } 42, p. 43-47.
\bibitem{RoVo} Y. A. Rozanov and V. A. Volkonskii (1959). Some limit
  theorems for random functions I. {\it Theory Probab. Appl.} 4, p. 178-197.
\bibitem{St} C. J. Stone (1966). On local and ratio limit theorems. {\it
  Proc. Fifth Berkeley Sympos. Math. Statist. Probab.}, Univ. Calif.,
  p. 217-224.
\bibitem{Ta} Y. Taga (1965). The optimal sampling procedure for estimating
  the mean of stationary Markov process. {\it Ann. Inst. Stat. Math.} 17,
  p. 105-112.
\bibitem{Utev0} S. A. Utev (1989). Sums of random variables with
  $\varphi$-mixing condition. In: {\it Asymptotical Analysis of Stochastic
    Processes Distributions}, Novosibirsk.
\bibitem{Utev} S. A. Utev (1990). Central limit theorem for dependent random
  variables. {\it Probab. Theory Math. Statist.} 2, p.
  519-528.
  \bibitem{Utevbis} S. A. Utev (1991). Sums of random variables with
  $\varphi$-mixing. {\it Siberian Advances in Mathematics} 1, 3,
  p.124-155.
\bibitem{With} C. S. Withers (1981). Central limit theorems for dependent
  variables. I. {\it Z. Wahrsch. Verw. Gebiete}  57, 4,
  p. 509-534. (Corrigendum in {\it Z. Wahrsch. Verw. Gebiete} 63, 4, (1983),
  p. 555).
\end{thebibliography}
\end{document}